\documentclass[11pt]{article}

\usepackage{amsfonts}
\usepackage{amssymb}
\usepackage{amsmath}
\usepackage{amsthm}
\def\negthickspace{\!\!\!}

\newcommand{\nicefrac}[2]
{\leavevmode \kern.1em\raise.5ex\hbox{\the\scriptfont0 #1}
             \kern-.1em/\kern-.15em\lower.25ex
             \hbox{\the\scriptfont0 #2}}

\newtheorem*{theorem}{Theorem}
\newtheorem*{proposition}{Proposition}

\newtheorem*{definition}{Definition}

\theoremstyle{definition}
\newtheorem*{remark}{Remark}
\newtheorem*{remarks}{Remarks}

\theoremstyle{definition}

\begin{document} 

\begin{center}
{\Large{\sc On critical normal sections}}\\[0.2cm]
{\Large{\sc for two-dimensional immersions in $\mathbb R^4$}}\\[0.2cm]
{\Large{\sc and a Riemann-Hilbert problem}}\\[1cm]
{\large Steffen Fr\"ohlich,\quad Frank M\"uller}\\[0.4cm]
{\small\bf Abstract}\\[0.4cm]
\begin{minipage}[c][2.5cm][l]{12cm}
{\small For orthonormal normal sections of two-dimensional immersions in $\mathbb R^4$ we define torsion coefficients and a functional for the total torsion. We discuss normal sections which are critical for this functional. In particular, a global estimate for the torsion coefficients of a critical normal section in terms of the curvature of the normal bundle is provided.}
\end{minipage}
\end{center}
{\small MCS 2000: 53A07, 53A10, 30G20}


\section{Introduction}
Consider a two-dimensional, conformally parametrized immersion
\begin{equation}\label{1.1}
  X=X(u,v)=\big(x^1(u,v),x^2(u,v),x^3(u,v),x^4(u,v)\big)\in C^4(B,\mathbb R^4)
\end{equation}
on the closed unit disc $B=\big\{(u,v)\in\mathbb R^2\,:\,u^2+v^2\le 1\big\}\subset\mathbb R^2,$ together with an orthogonal moving $4$-frame
\begin{equation}\label{1.2}
  \{X_u,X_v,N_1,N_2\},
\end{equation}
which consists of the orthogonal tangent vectors $X_u,$ $X_v,$ and orthogonal unit normal vectors $N_1,$ $N_2\in C^3(B,\mathbb R^4)$:
\begin{equation}\label{1.3}
\begin{array}{l}
  X_u\cdot X_u^t=:g_{11}=W=g_{22}:=X_v\cdot X_v^t\,,\quad
  g_{12}:=X_u\cdot X_v^t=0, \\[0.2cm]
  X_u\cdot N_\sigma^t=0=X_v\cdot N_\sigma^t\quad\mbox{for}\ \sigma=1,2, \\[0.2cm]
  |N_1|=1=|N_2|,\quad
  N_1\cdot N_2^t=0.
\end{array}
\end{equation}
Here $W$ denotes the area element of $X,$ and $X^t$ means the transposed vector of $X.$\\[1ex]
Finally, we set $\mathring{B}=\big\{(u,v)\in\mathbb R^2\,:\,u^2+v^2<1\big\}$ for the open unit disc and $\partial B=\big\{(u,v)\in\mathbb R^2\,:\,u^2+v^2=1\big\}$ for its boundary.
\begin{remarks}\quad
\begin{itemize}
\item[1.]
Note the relation $W>0$ in $B.$
\vspace*{-1.2ex}
\item[2.]
For the introduction of conformal parameters into a Riemannian metric we refer to \cite{Sauvigny_01}.
\end{itemize}
\end{remarks}
\noindent
With the present notes we want to draw the reader's attention to a definition of torsion coefficients $T_{\sigma,i}^\vartheta$ for an orthonormal normal section $\{N_1,N_2\}$, which is deduced from the theory of space curves. In addition, we introduce an associated functional of total torsion ${\mathcal T}_X(N_1,N_2)$ and study its critical points.
\goodbreak\noindent
The paper is organized as follows:
\vspace*{-0.7ex}
\begin{itemize}
\item
In Chapter 2, we define torsion coefficients of orthonormal normal sections.
\vspace*{-1.2ex}
\item
In Chapter 3, we introduce the concept of the total torsion of an orthonormal normal section. We provide conditions for a normal section to be critical and optimal for the total torsion.
\vspace*{-1.2ex}
\item
Chapter 4 contains some aspects about generalized analytic functions and Riemann-Hilbert problems. These results are used to prove a global pointwise estimate for the torsion coefficients of critical normal sections. 
\vspace*{-1.2ex}
\item
Finally, an example of a critical normal section for holomorphic graphs $(w,\Phi(w))$ will be discussed in Chapter 5.
\end{itemize}


\section{Torsion coefficients and curvature of the normal bundle}
According to the classical theory of curves in $\mathbb R^3,$ we introduce \emph{torsion coefficients} as follows:
\begin{definition}
For an orthonormal normal section $\{N_1,N_2\}$ we define
\begin{equation}\label{2.1}
  T_{\sigma,i}^\vartheta:=N_{\sigma,u^i}\cdot N_\vartheta^t\,,\quad i=1,2,\ \sigma,\vartheta=1,2,
\end{equation}
setting $u^1\equiv u$ and $u^2\equiv v.$
\end{definition}
\begin{remarks}\quad
\begin{itemize}
\item[1.]
Obviously, $T_{\sigma,i}^\vartheta=-T_{\vartheta,i}^\sigma$ holds for any $i=1,2,$ $\sigma,\vartheta=1,2,$ and consequently $T_{\sigma,i}^\sigma\equiv 0.$
\vspace*{-1.2ex}
\item[2.]
The $T_{\sigma,i}^\vartheta$ are exactly the coefficients of the normal connection. In their terms one defines the coefficients of the curvature tensor ${\mathfrak S}$ of the normal bundle (summation convention!)
\begin{equation}\label{2.2}
  S_{\sigma,ij}^\vartheta
  :=T_{\sigma,i,u^j}^\vartheta-T_{\sigma,j,u^i}^\vartheta
   +T_{\sigma,i}^\omega T_{\omega,j}^\vartheta-T_{\sigma,j}^\omega T_{\omega,i}^\vartheta\,,\quad
  i,j=1,2,\ \sigma,\vartheta=1,2.
\end{equation}
In contrast to the case of codimension $n\ge 3,$ the quadratical terms in (\ref{2.2}) vanish in $\mathbb R^4,$ and ${\mathfrak S}$ consists essentially of the single term
\begin{equation}\label{2.3}
  S:=S_{1,12}^2
  =T_{1,1,v}^2-T_{1,2,u}^2
  =\mbox{div}\,(-T_{1,2}^2,T_{1,1}^2).\vspace*{-1.2ex}
\end{equation}
This is the reason why we concentrate on immersions in $\mathbb R^4.$
\vspace*{-1.2ex}
\item[3.]
Note that $S$ does not depend on the choice of the orthonormal section $\{N_1,N_2\}$, compare Subsection\,3.1.
\vspace*{-1.2ex}
\item[4.]
Distinguish our definition from that of the {\it normal torsion} of a surface (see i.e.~\cite{Faessler_01}): It can be defined as the torsion of the one-dimensional normal section (as a curve on the surface), which arises from a suitable intersection of $X$ with a three-dimensional hyperplane.
\end{itemize}
\end{remarks}
\noindent


\section{Total torsion and optimal normal sections}
\setcounter{equation}{0}
To an orthonormal section $\{N_1,N_2\}$ of the normal bundle we assign the {\it total torsion}
\begin{equation}\label{3.1}
  {\mathcal T}_X(N_1,N_2)
  :=\sum_{\sigma,\vartheta=1}^2\,
   \int\hspace*{-0.25cm}\int\limits_{\hspace{-0.3cm}B}
   g^{ij}\,T_{\sigma,i}^\vartheta T_{\sigma,j}^\vartheta\,W\,dudv
  =2\int\hspace*{-0.25cm}\int\limits_{\hspace{-0.3cm}B}
    \Big\{(T_{1,1}^2)^2+(T_{1,2}^2)^2\Big\}\,dudv,
\end{equation}
where $g_{ij}g^{jk}=\delta_i^k,$ and $\delta_i^k$ is the Kronecker symbol (see the conformality relations in (\ref{1.3})).

\subsection{Critical orthonormal normal sections}
The total torsion depends on the chosen orthonormal section $\{N_1,N_2\},$ and it can be controlled by means of a rotation angle $\varphi=\varphi(u,v),$ depending smoothly on $(u,v)\in B.$ Indeed, starting with the section $\{N_1,N_2\},$ we write 
\begin{equation}\label{3.2}
  \widetilde N_1=\cos\varphi\,N_1+\sin\varphi\,N_2\,,\quad
  \widetilde N_2=-\sin\varphi\,N_1+\cos\varphi N_2\,
\end{equation}
for the rotated normal section $\{\widetilde N_1,\widetilde N_2\}$. Then, the new torsion coefficients are given by
\begin{equation}\label{3.3}
  \widetilde T_{1,1}^2=T_{1,1}^2+\varphi_u\,,\quad
  \widetilde T_{1,2}^2=T_{1,2}^2+\varphi_v\,.  
\end{equation}
Due to (\ref{3.1}), the difference between new and old total torsion now computes to
\begin{equation}\label{3.4}
\begin{array}{l}
  \displaystyle
  {\mathcal T}_X(\widetilde N_1,\widetilde N_2)-{\mathcal T}_X(N_1,N_2)
  \,=\,2\int\hspace{-0.25cm}\int\limits_{\hspace{-0.3cm}B}|\nabla\varphi|^2\,dudv
       +4\int\hspace{-0.25cm}\int\limits_{\hspace{-0.3cm}B}
         (T_{1,1}^2\varphi_u+T_{1,2}^2\varphi_v)\,dudv \\[0.8cm]
  \hspace*{0.6cm}\displaystyle
  =\,2\int\hspace{-0.25cm}\int\limits_{\hspace{-0.3cm}B}|\nabla\varphi|^2\,dudv
     +4\int\limits_{\partial B}(T_{1,1}^2,T_{1,2}^2)\cdot\nu^t\,\varphi\,ds
     -4\int\hspace{-0.25cm}\int\limits_{\hspace{-0.3cm}B}
       \mbox{div}\,(T_{1,1}^2,T_{1,2}^2)\varphi\,dudv.
\end{array}
\end{equation}
In general, the right hand side does not vanish.
\begin{proposition}
Let $\{N_1,N_2\}$ be critical for ${\mathcal T}_X.$ Then the torsion coefficients satisfy
\begin{equation}\label{3.5}
  \mbox{\rm div}\,(T_{1,1}^2,T_{1,2}^2)=0\quad\mbox{in}\ B,\quad
  (T_{1,1}^2,T_{1,2}^2)\cdot\nu^t=0\quad\mbox{on}\ \partial B.
\end{equation}
\end{proposition}

\subsection{Construction of critical orthonormal normal sections}
How can we construct a critical section $\{N_1,N_2\}$ from a given section $\{\widetilde N_1,\widetilde N_2\}?$\\[1ex]
If $\{N_1,N_2\}$ is critical, then we have
\begin{equation}\label{3.8}
\begin{array}{l}
  0\,=\,\mbox{div}\,(T_{1,1}^2,T_{1,2}^2)
   \,=\,\mbox{div}\,(\widetilde T_{1,1}^2-\varphi_u,\widetilde T_{1,2}^2-\varphi_v)\quad\mbox{in}\ B, \\[0.2cm]
  0\,=\,(T_{1,1}^2,T_{1,2}^2)\cdot\nu^t
   \,=\,(\widetilde T_{1,1}^2-\varphi_u,\widetilde T_{1,2}^2-\varphi_v)\cdot\nu^t\quad\mbox{on}\ \partial B,
\end{array}
\end{equation}
by virtue of (\ref{3.3}), (\ref{3.5}). This implies our next result:
\begin{proposition}
The given section $\{\widetilde N_1,\widetilde N_2\}$ transforms into a critical section by means of (\ref{3.2}), iff
\begin{equation}\label{3.9}
\begin{array}{l}
  \Delta\varphi=\mbox{\rm div}\,(\widetilde T_{1,1}^2,\widetilde T_{1,2}^2)\quad\mbox{in}\ B, \\[0.4cm]
  \displaystyle
  \frac{\partial\varphi}{\partial\nu}=(\widetilde T_{1,1}^2,\widetilde T_{1,2}^2)\cdot\nu^t
  \quad\mbox{on}\ \partial B
\end{array}
\end{equation}
holds for the rotation angle $\varphi=\varphi(u,v)$.
\end{proposition}
\begin{remark}
It is well known that the solvability of the Neumann problem
\begin{equation}\label{3.10}
    \Delta\varphi=f\quad\mbox{in}\ B,\quad
    \frac{\partial\varphi}{\partial\nu}=g\quad\mbox{on}\ \partial B
\end{equation}
depends on the integrability condition
\begin{equation}\label{3.11}
  \int\hspace*{-0.25cm}\int\limits_{\hspace*{-0.3cm}B}f\,dudv
  =\int\limits_{\partial B}g\,ds,
\end{equation}
which is fulfilled in our proposition.
\end{remark}
\subsection{Minimality of critical orthonormal normal sections}
\noindent
Let $\{N_1,N_2\}$ be a critical section. Then, we conclude
\begin{equation}\label{3.12}
  {\mathcal T}_X(\widetilde N_1,\widetilde N_2)
  ={\mathcal T}_X(N_1,N_2)
   +2\int\hspace{-0.25cm}\int\limits_{\hspace{-0.3cm}B}
     |\nabla\varphi|^2\,dudv,
\end{equation}
taking (\ref{3.4}) and (\ref{3.5}) into account. This proves the following
\begin{proposition}
A critical orthonormal normal section $\{N_1,N_2\}$ minimizes the total torsion, i.e.~we have
\begin{equation}\label{3.13}
  {\mathcal T}_X(N_1,N_2)\le{\mathcal T}_X(\widetilde N_1,\widetilde N_2)
\end{equation}
for all smooth orthonormal normal sections $\{\widetilde N_1,\widetilde N_2\}.$ The equality occurs iff $\varphi\equiv\mbox{const}.$
\end{proposition}
\subsection{Flat normal bundles}
For a critical normal section, the vector-field $(-T_{1,2}^2,T_{1,1}^2)$ is parallel to $\nu$ along $\partial B$. Applying the Gaussian integral theorem to (\ref{2.3}), we infer
\begin{equation}\label{3.6}
  \int\hspace*{-0.25cm}\int\limits_{\hspace*{-0.3cm}B}S\,dudv=\int\limits_{\partial B}(-T_{1,2}^2,T_{1,1}^2)\cdot\nu^t\,ds
  =\pm\int\limits_{\partial B}\sqrt{(T_{1,1}^2)^2+(T_{1,2}^2)^2}\,ds.
\end{equation}
In particular, if $S\equiv 0,$ that is, the normal bundle is flat, then we find
\begin{equation}\label{3.7}
  T_{\sigma,i}^\vartheta\equiv 0\quad\mbox{on}\ \partial B
\end{equation}
for $i=1,2$ and $\sigma,\vartheta=1,2.$\\[1ex]
Differentiating (\ref{2.3}) and (\ref{3.5}), we further obtain
\begin{equation}\label{3.14}
  \Delta T_{1,1}^2=\frac{\partial}{\partial v}\,S=0\,,\quad
  \Delta T_{1,2}^2=-\,\frac{\partial}{\partial u}\,S=0
  \quad\mbox{in}\ B
\end{equation}
for flat normal bundles. Therefore,
\begin{equation}\label{3.15}
  T_{\sigma,i}^\vartheta\equiv 0\quad\mbox{in}\ B\quad(i=1,2,\ \ \sigma,\vartheta=1,2)
\end{equation}
follows by the maximum principle.
\begin{remark}
Immersions of prescribed mean curvature with flat normal bundles are extensively studied in the literature; see e.g. \cite{Smoczyk_Wang_Xin_01}, \cite{Froehlich_Winklmann_01} for higher dimensional surfaces. Special results for two-dimensional immersions without curvature conditions on the normal bundle can be found in \cite{Froehlich_04}.
\end{remark}
\noindent
In the following, we investigate the inhomogeneous case of non-flat normal bundles to extend the relation (\ref{3.15}) appropriately.
\section{Estimates for the torsion coefficients}
\setcounter{equation}{0}
\subsection{A Riemann-Hilbert problem}
Once again, let us consider (\ref{2.3}) and (\ref{3.5}) for critical sections:
\begin{equation}\label{4.1}
  \frac{\partial}{\partial u}\,T_{1,1}^2+\frac{\partial}{\partial v}\,T_{1,2}^2=0,\quad
  \frac{\partial}{\partial v}\,T_{1,1}^2-\frac{\partial}{\partial u}\,T_{1,2}^2=S
  \quad\mbox{in}\ B.
\end{equation}
The {\it complex-valued torsion} $\Psi=T_{1,1}^2-iT_{1,2}^2$ solves the non-homogeneous Cauchy-Rie\-mann equation
\begin{equation}\label{4.2}
  \frac{\partial}{\partial\overline w}\,\Psi(w)=\Psi_{\overline w}(w)
  :=\frac{1}{2}\,(\Psi_u+i\Psi_v)
  =\frac{i}{2}\,S,\quad w=u+iv\in\mathring{B}.
\end{equation}
In addition, we write the boundary condition in (\ref{3.5}) as
\begin{equation}\label{4.3}
  \mbox{Re}\big[w\Psi(w)\big]=0,\quad w\in\partial B.
\end{equation}
The relations (\ref{4.2}) and (\ref{4.3}) form a linear Riemann-Hilbert problem for $\Psi.$
\begin{proposition}
The problem (\ref{4.2})-(\ref{4.3}) possesses at most one solution $\Psi\in C^1(\mathring{B})\cap C^0(B)$.
\end{proposition}
\begin{proof}
Let $\Psi_1,\Psi_2$ be two such solutions. Then we set $\Phi(w):=w[\Psi_1(w)-\Psi_2(w)]$ and note
\begin{equation}\label{4.4}
  \Phi_{\overline w}=0\quad\mbox{in}\ \mathring{B},\qquad \mbox{Re}\,\Phi=0\quad\mbox{on}\ \partial B.
\end{equation}
Consequently, $\Phi\equiv ic$ holds true in $B$ with some constant $c\in\mathbb R$, and the continuity of $\Psi_1,\Psi_2$ implies $c=0$.
\end{proof}
\subsection{Some facts about generalized analytic functions}
As general references for this subsection we name \cite{Vekua}, \cite{Sauvigny_02}.\\[1ex]
For arbitrary $f\in C^1(B,\mathbb C)$ we define 
\begin{equation}\label{4.5}
  T_B[f](w)
  :=-\frac1\pi
     \int\hspace*{-0.25cm}\int\limits_{\hspace*{-0.3cm}B}
     \frac{f(\zeta)}{\zeta-w}\,d\xi d\eta,\quad w\in\mathbb C,
\end{equation}
using the notation $\zeta=\xi+i\eta$. Then, there hold $g:=T_B[f]\in C^1(\mathbb C\setminus\partial B)\cap C^0(\mathbb C)$ as well as 
\begin{equation}\label{4.6}
  \frac\partial{\partial\overline w}\,T_B[f](w)=\left\{
  \begin{array}{ll}
  f(w),& w\in\mathring{B}\\[1ex]
  0,& w\in\mathbb C\setminus B
  \end{array}\right.,
\end{equation}
cf.~\cite{Vekua} Kapitel I, \S5. Next, we set
\begin{equation}\label{4.7}
\begin{array}{rcl}
  P_B[f](w)\negthickspace
  & := & \negthickspace\displaystyle
         -\,\frac1\pi
            \int\hspace*{-0.25cm}\int\limits_{\hspace*{-0.3cm}B}
            \left\{
              \frac{f(\zeta)}{\zeta-w}+\frac{\overline\zeta\,
              \overline{f(\zeta)}}{1-w\overline\zeta}
            \right\}d\xi\,d\eta\\[4ex]
  &  = & \negthickspace\displaystyle
         T_B[f](w)+\frac1w\,\overline{T_B[wf]\Big(\frac1{\overline w}\Big)}\,.
\end{array}
\end{equation}
We obtain $h:=P_B[f]\in C^1(\mathring{B})\cap C^0(B)$, and (\ref{4.6}) yields
\begin{equation}\label{4.8}
  \frac\partial{\partial\overline w}\,P_B[f](w)=f(w),\quad w\in\mathring{B}.
\end{equation}
Finally, we note the relation
\begin{equation}\label{4.9}
  P_B[f](w)=T_{\mathbb C}[f_*](w),\quad w\in B.
\end{equation}
Here, $T_{\mathbb C}$ is defined as $T_B$ but with integration over $\mathbb C$, and we have abbreviated
\begin{equation}\label{4.10}
  f_*(w):=\left\{
  \begin{array}{ll}
  f(w),& w\in B\\[1ex]
  \displaystyle\frac1{|w|^4}\,\overline{f\Big(\frac 1{\overline w}\Big)}, & w\in\mathbb C\setminus B
  \end{array}\right..
\end{equation}
Observe that $f_*(w)$ is not continuous in $\mathbb C$, but it belongs to the class $L_{p,2}(\mathbb C)$ for any $p\in[1,+\infty]$, that means, $f_*(w)$ as well as $|w|^{-2}f_*(\frac 1w)$ belong to $L_p(B)$; compare \cite{Vekua} p.\,12. Consequently, Satz\,1.24 in \cite{Vekua} yields the following
\begin{proposition}
With the definitions above, we have the uniform estimate
\begin{equation}\label{4.11}
  \big|P_B[f](w)\big|=\big|T_{\mathbb C}[f_*](w)\big|\le c(p)\|f\|_{L_p(B)},\quad w\in B,
\end{equation}
where $p\in(2,+\infty],$ and $c(p)$ is a positive constant dependent on $p$.
\end{proposition}

\subsection{A global pointwise estimate for the torsion coefficients}
\begin{theorem}
Consider a conformally parametrized immersion $X\in C^4(B,\mathbb R^4)$ and write
\begin{equation}\label{4.12}
  s_p:=\|S\|_{L_p(B)},\quad p\in(2,+\infty].
\end{equation}
Then, the complex-valued torsion $\Psi=T_{1,1}^2-iT_{1,2}^2$ of a critical orthonormal section $\{N_1,N_2\}$ satisfies
\begin{equation}\label{4.13}
  |\Psi(w)|\le c(p)s_p\quad\mbox{for all}\ w\in B,
\end{equation}
with some positive constant $c(p)$.
\end{theorem}
\begin{remark}\quad
\begin{itemize}
\item[1.]
For a flat normal bundle, i.e. $s_p=0$, we recover (\ref{3.15}).
\vspace*{-1.2ex}
\item[2.]
The general estimate (\ref{4.13}) shall be useful, e.g., for proving curvature estimates for immersions with non-flat normal bundle.
\end{itemize}
\end{remark}
\begin{proof}[Proof of the theorem]
Let us write $f:=\frac i2S\in C^1(B)$. We claim that $\Psi$ possesses the integral representation 
\begin{equation}\label{4.15}
  \Psi(w)
  =P_B[f](w)
  =-\frac{1}{\pi}\,
    \int\hspace{-0.25cm}\int\limits\limits_{\hspace{-0.3cm}B}
    \left\{
      \frac{f(\zeta)}{\zeta-w}+\frac{\overline\zeta\,\overline{f(\zeta)}}{1-w\overline\zeta}	
    \right\}d\xi\,d\eta,
  \quad w\in B.
\end{equation}
Then, (\ref{4.13}) follows at once from the proposition in Subsection\,4.2.\\[1ex]
An elementary calculation proves
\begin{equation}\label{4.16}
  wP_B[f](w)
  =\frac{1}{\pi}\,
   \int\hspace{-0.25cm}\int\limits_{\hspace{-0.3cm}B}
   f(\zeta)\,d\xi\,d\eta+T_B[wf](w)-\overline{T_B[wf]\Big(\frac 1{\overline w}\Big)}\ .
\end{equation}
Taking $f=\frac i2S$ into account, we infer
\begin{equation}\label{4.17}
  \mbox{Re}\big\{wP_B[f](w)\big\}=0,\quad w\in\partial B.
\end{equation}
Consequently -- remember (\ref{4.8}) --, $P_B[f](w)$ solves the Riemann-Hilbert problem (\ref{4.2})-(\ref{4.3}). Now the uniqueness result of the proposition in Subsection\,4.1 yields the identity (\ref{4.15}).
\end{proof}
\begin{remark}
We point out that the representation (\ref{4.15}) relies crucially on the fact that the right hand side $f=\frac i2S$ in (\ref{4.2}) is purely imaginary.  In general, a Riemann-Hilbert problem as in (\ref{4.2})-(\ref{4.3}) is solvable, iff the integral of the right hand side $f$ over $B$ has vanishing real part, cf.~(\ref{4.16}). For details we refer to \cite{Vekua} Kapitel IV, \S7.
\end{remark}


\section{Example: Holomorphic graphs on $B$}
\setcounter{equation}{0}
Let us consider graphs $X(w)=(w,\Phi(w)),$ $w=u+iv\in B.$ If $\Phi(w)=\varphi(w)+i\psi(w)$ is holomorphic on $B$, then the vectors 
\begin{equation}\label{5.1}
  N_1=\frac{1}{\sqrt{W}}\,(-\varphi_u,-\varphi_v,1,0),\quad
  N_2=\frac{1}{\sqrt{W}}\,(-\psi_u,-\psi_v,0,1)
\end{equation}
form an orthonormal normal section, where $W=1+|\nabla\varphi|^2=1+|\Phi'|^2$ is the area element.
\begin{remark}
Due to $\varphi_u=\psi_v$, $\varphi_v=-\psi_u$ and thus $\Delta\varphi=\Delta\psi=0,$ the immersion $X$ represents a conformally parametrized minimal graph.
\end{remark}
\noindent
For the torsion coefficients we compute
\begin{equation}\label{5.2}
  T_{1,1}^2
  =\frac{1}{W}\,(-\varphi_{uu}\varphi_v+\varphi_{uv}\varphi_u)
  =\frac1{2W}\frac\partial{\partial v}(|\nabla\varphi|^2)\,,\quad
  T_{1,2}^2=-\frac1{2W}\frac\partial{\partial u}(|\nabla\varphi|^2)\,,
\end{equation}
Consequently, the relation
\begin{equation}\label{5.3}
  \mbox{div}\,(T_{1,1}^2,T_{1,2}^2)=0\quad\mbox{in}\ B
\end{equation}
is satisfied. In order to check the boundary condition in (\ref{3.5}), we introduce polar coordinates $u=r\cos\alpha,$ $v=r\sin\alpha$ and note $\frac1r\frac\partial{\partial\alpha}=u\frac\partial{\partial v}-v\frac\partial{\partial u}$. According to (\ref{5.2}), we then obtain
\begin{equation}\label{5.4}
  (T_{1,1}^2,T_{1,2}^2)\cdot\nu^t=\frac1{2W}\Big(u\frac\partial{\partial v}-v\frac\partial{\partial u}\Big)(|\nabla\varphi|^2)
  =\frac1{2W}\frac\partial{\partial\alpha}(|\Phi'|^2)\quad\mbox{on}\ \partial B.
\end{equation}  
\begin{proposition}
Consider the graph $(w,\Phi(w)),$ $w\in B,$ with a holomorphic function $\Phi(w)=\varphi(w)+i\psi(w)$. Then the normal section $\{N_1,N_2\}$ defined in (\ref{5.1}) is critical, that is, it satisfies (\ref{3.5}), iff $|\Phi'|$ is constant on $\partial B.$
\end{proposition}
\begin{remark}
As an example, we mention the graph $X(w)=(w,w^n)$, $w\in B$, for arbitrary $n\in\mathbb N$.
\end{remark}


\vspace*{0.8cm}
\noindent
Steffen Fr\"ohlich\\
Freie Universit\"at Berlin\\ 
Fachbereich Mathematik und Informatik\\ 
Arnimallee 2-6\\
D-14195 Berlin\\
Germany\\[0.2cm]
e-mail: sfroehli@mi.fu-berlin.de\\[1cm]
Frank M\"uller\\
Brandenburgische Technische Universit\"at Cottbus\\
Mathematisches Institut\\
Konrad-Zuse-Stra{\ss}e 1\\
D-03044 Cottbus\\
Germany\\[0.2cm]
e-mail: mueller@math.tu-cottbus.de

\end{document}